\documentclass[11pt, leqno]{amsart}
\usepackage{amsmath, amssymb, latexsym}
\usepackage{cases}
\usepackage{mathrsfs}
\usepackage{enumerate}
\usepackage{color}
\usepackage[colorlinks=true, pdfstartview=FitV, linkcolor=blue,%
citecolor=blue, urlcolor=blue]{hyperref}
\usepackage[all, knot]{xy}

\xyoption{arc}

\numberwithin{equation}{section}

\newtheorem{theorem}{Theorem}[section]
\newtheorem{corollary}[theorem]{Corollary}
\newtheorem{lemma}[theorem]{Lemma}
\newtheorem{proposition}[theorem]{Proposition}

\newtheorem{example}[theorem]{Example}

\theoremstyle{definition}
\newtheorem{definition}[theorem]{Definition}

\newcommand{\Q}{\mathbb{Q}}

\newcommand{\Z}{\mathbb{Z}}
\newcommand{\N}{\mathbb{N}}
\newcommand{\h}{\mathfrak{h}}
\newcommand{\g}{\mathfrak{g}}

\title[Braid group actions of quantum Borcherds-Bozec algebras]
{ Braid group actions of quantum Borcherds-Bozec algebras}

\author[Zhaobing Fan]{Zhaobing Fan}
\address{Harbin Engineering University,
Harbin, China}
\email{fanzhaobing@hrbeu.edu.cn}
\thanks{ }

\author[Bolun Tong]{Bolun Tong}
\address{Harbin Engineering University,
Harbin, China}
\email{tbl\_2019@hrbeu.edu.cn}

\address{}


\begin{document}

\begin{abstract}

In this paper, we construct the Lusztig symmetries for quantum Borcherds-Bozec algebra $U_q(\g)$ and its weight module $M\in \mathcal O$, on which the generators with real indices of $U_q(\g)$  act nilpotently. We show that these symmetries satisfy the defining relations of the braid group, associated to the Weyl group $W$ of $U_q(\g)$, which gives a braid group action.

\end{abstract}

\maketitle

\section*{Introduction}

\vskip 2mm
 Lusztig considered in \cite{Lus90} certain perverse sheaves on representation varities of a quiver of type A,D,E, and gave a geometric approach to the half parts of corresponding quantum groups. 
The canonical basis theory arose in this setting, which is given by simple perverse sheaves. Moreover, 
the algebraic and coalgebraic structures of the quantum groups are given by induction functor and restriction functor, respectively.
 These results were generalized by  
Lusztig in \cite{Lus91} to Kac-Moody cases. 
Later on, Lusztig consider more general cases, namely arbitrary quivers, possibly carrying loops in \cite{Lus93}.
The obtained algebra is denoted by $U^{-}$.
In \cite{Lus93}, Lusztig proposed a question that if $U^{-}$ is generated by the elementary simple perverse sheaves $F_i^{(n)}$ with all vertices $i$ and $n\in \N$ as an algebra.
 The question is answered by himself in the case of the quiver with one vertex and multiple loops, by a quadratic form criterion for a monomial to be {\it tight} or {\it semi-tight}. 
Based on  Kang and Schiffmann's  work for quantum generalied Kac-Moody algebras (cf. \cite{KS06}), 
 Li and Lin \cite{Lin} answered the question when the quiver has at least two loops on each imaginary vertex.
   In \cite{Bozec2014b}, Bozec solved the Lusztig's question completely. As a bialgebra, the resulting $U^-$ is so called {\it quantum Borcherds-Bozec algebra}, which is the main object we studied in the current paper.

\vskip 2mm
On algebraic side, the quantum Borcherds-Bozec algebras can be treated as a further generalization of quantum generalized Kac-Moody algebras \cite{Kang95}. More precisely, a  quantum Borcherds-Bozec algebra $U_q(\g)$ has infinitely many generators $e_{il}, f_{il}$ $(l \in \Z_{>0})$ for each imaginary index $i\in I^{\text{im}}$, and their degrees are $l$ multiples of $\alpha_i$ and $-\alpha_i$, respectively. The commutation relations between these generators are rather complicated and are higher order in some sense (cf. \cite{FKKT}).
Thanks to Bozec, there exists a set of {\it primitive generators} $\mathtt a_{il}, \mathtt b_{il}$ $((i,l)\in I^\infty)$  with better properties and simpler commutation relations. Using these generators, Bozec constructed the Kashiwara operators. He then developed the crystal basis theory for quantum Borcherds-Bozec algebras and their irreducible highest weight modules \cite{Bozec2014c}. In  \cite{FKKT21}, the authors and Kang and Kim constructed the Global bases for the quantum Borcherds-Bozec algebras.

\vskip 2mm

Lusztig defined in \cite{Lusztig} the {\it symmetries} $T'_{i,e}, T''_{i,e}$ ($e=\pm 1$) for the quantum group $U$ of Kac-Moody type and its integrable weight modules. He proved that these automorphisms $T'_{i,e}, T''_{i,e}$  satisfy the braid group relations both on $U$ and integrable $U$-modules.
As an important application, he gave  some linearly independent subsets of $U$ associated to reduced expressions of elements in the Weyl group $W$. In the case of finite types, if one choose the longest element $\omega_0\in W$, the  independent set considered form the PBW-basis. He also shown how to extend  as much as possible this construction to arbitrary Cartan data, especially for affine cases.

\vskip 2mm

In this paper, we shall define the Lusztig symmetries for quantum Borcherds-Bozec algebra $U_q(\g)$. Note that, the Weyl group of $U_q(\g)$ is generated by the reflections associated to real indexes, it is not enough to construct the PBW-type basis as Lusztig did. But since the primitive  generators $\mathtt a_{il}$ and $\mathtt b_{il}$ satisfy the `higer order' quantum Serre relations for $i\in I^{\text{im}}$, we could define the Lusztig symmetries in a natural way.
In the case where $I$ consists of exact two real indexes $i\neq j$ and $a_{ij}a_{ji}$ is finite, Lusztig proved that there is a braid group action on the integral modules (cf. \cite{Lusztig}) through the symmetries. In Theorem \ref{Ba}, we prove that it can be generalized to an arbitrary Borcherds-Cartan datum by using our constructed symmetries, which give the braid group actions on $U_q(\g)$ and $U_q(\g)$-module $M$ in a certain category $\mathcal O$.

\vskip 2mm
In Lusztig's construction for PBW-basis, a crucial result is that the inner product of $U^+$ is $T''_{i,1}$-invariant on the subalgebra $U^+[i]$, which is generated by a set of elements of the forms similar to the Serre-type relations  but allows for smaller degrees. 
To verify this in quantum Borcherds-Bozec algebra case (Theorem \ref{remain}), we follow the framework given in \cite[Chapter 8A]{Jantzen} rather that Section 38.2 in \cite{Lusztig}, since we have a more general setting for the values of our bilinear form. We make a notice here that, as a consequence of Theorem \ref{remain}, one could get a lot of linearly independent subsets of $U^+_q(\g)$ as in \cite[38.2.2]{Lusztig} by a similar argument.

\vskip 2mm

This paper is organized as follows.
In Section 1, we review the definition of quantum Borcherds-Bozec algebras and the notion of the primitive generators. In Section 2, we define the Lusztig symmetries $L'_{i,e}, L''_{i,e}$ ($e=\pm 1$) for these algebras and their weight module $M\in \mathcal O$, and prove the braid group actions on them. In Section 3, we investigate the relations between the symmetries and the bilinear form $\{ \ , \ \}$ on $U^+_q(\g)\times U^-_q(\g)$.

\vskip 2mm

\noindent\textbf{Acknowledgements.}

Z. Fan was partially supported by the NSF of China grant 11671108, the NSF of Heilongjiang Province grant JQ2020A001, and the Fundamental Research Funds for the central universities. We would like to express our sincere gratitude to Professor Seok-Jin Kang for his helpful discussions.

\vskip 2mm

\section{Quantum Borcherds-Bozec algebras}
\vskip 2mm

Let $I$ be a finite or countably infinite index set. An integer-valued matrix $A=(a_{ij})_{i,j\in I}$ is called an {\it even symmetrizable Borcherds-Cartan matrix} if it satisfies the following conditions:
\begin{itemize}
\item[(i)] $a_{ii}=2,0,-2,-4,\dots$,
\item[(ii)] $a_{ij}\in \Z_{\leq 0}$ for $i\neq j$,
\item[(iii)] there is a diagonal matrix $D=\text{diag}(s_i\in \Z_{>0}\mid i\in I)$ such that $DA$ is symmetric.
\end{itemize}

\vskip 3mm

Let $I^{\text{re}}=\{i\in I\mid a_{ii}=2\}$ be the set of {\it real indices}. Let $I^{\text{im}}=\{i\in I\mid a_{ii}\leq 0\}$ and $I^{\text{iso}}=\{i\in I\mid a_{ii}= 0\}$ be the set of {\it imaginary indices} and {\it isotropic indices}, respectively.

\vskip 3mm
A {\it Borcherds-Cartan datum} consists of
\begin{itemize}
\item[(a)] an even symmetrizable Borcherds-Cartan matrix $A=(a_{ij})_{i,j\in I}$,
\item[(b)] a free abelian group $P^\vee=(\bigoplus_{i\in I}\Z h_i)\oplus (\bigoplus_{i\in I}\Z d_i)$, the {\it dual weight lattice},
\item[(c)] $\h=\Q\otimes_{\Z}P^\vee$, the {\it Cartan subalgebra},
\item[(d)] $P=\{\lambda \in \h^*\mid \lambda(P^\vee)\subseteq \Z\}$, the {\it weight lattice},
\item[(e)] $\Pi^\vee=\{h_i\in P^\vee\mid i\in I\}$, the set of {\it simple coroots},
\item[(f)] $\Pi=\{\alpha_i\in P \mid i\in I\}$, the set of {\it simple roots}, which is linearly independent over $\Q$ and satisfies
$$\alpha_j(h_i)=a_{ij},\ \alpha_j(d_i)=\delta_{ij}\ \ \text{for all}\ i,j \in I,$$
\item[(g)] for each $i\in I$, there is an element $\Lambda_i\in P$,  called the {\it fundamental weight}, defined by
$$\Lambda_i(h_j)=\delta_{ij},\ \Lambda_i(d_j)=0 \ \ \text{for all} \ i,j\in I.$$
\end{itemize}

\vskip 3mm

We denote by $P^+$ the set $\{\lambda\in P\mid \lambda(h_i)\geq 0 \ \text{for all}\ i\in I\}$ of {\it dominant integral weights}. The free abelian group $Q=\bigoplus_{i\in I}\Z\alpha_i$ is called the {\it root lattice}. Set $Q_+=\sum_{i\in I}\Z_{\geq 0}\alpha_i$ and $Q_-=-Q_+$. For $\beta=\sum_{i\in I} k_i\alpha_i\in Q_+$, we define its {\it hight} to be $|\beta|=\sum_{i\in I} k_i$.

\vskip 3mm

There is a non-degenerate symmetric bilinear form $( \ , \ )$ on $\h^*$ satisfying
$$(\alpha_i, \lambda)=s_i\lambda(h_i), \ (\Lambda_i,\lambda)=s_i\lambda(d_i) \ \ \text{for any} \ \lambda \in \h^* \ \text{and}\ i\in I,$$
and therefore we have
$$(\alpha_i,\alpha_j)=s_ia_{ij}=s_ja_{ji}\ \ \text{for all} \ i,j \in I.$$

\vskip 3mm
For $i\in I^{\text{re}}$, we define the {\it simple reflection} $r_i\in GL(\h^*)$ by
$$r_i(\lambda)=\lambda-\lambda(h_i)\alpha_i \  \ \text{for} \ \lambda \in \h^*.$$
The subgroup $W$ of $GL(\h^*)$ generated by $r_i$ $(i\in I^{\text{re}})$ is called the {\it Weyl group} of the Borcherds-Cartan datum. Note that the symmetric bilinear form $( \ , \ )$ is $W$-invariant.

\vskip 3mm

$W$ is a coxeter group generated by $r_i$ $(i\in I^{\text{re}})$ with defining relations $r_i^2=1$ $(i\in I^{\text{re}})$ and $(r_ir_j)^{m_{ij}}=1$ $(i\neq j)$, where $m_{ij}$ is the order of $r_ir_j$ and is related to $a_{ij}$ as follows:

\begin{figure}[h]

        \centering

       \begin{tabular}{|c|ccccc|}

         \hline


         $a_{ij}a_{ji}$ & $0$ & $1$ & $2$ & $3$ & $\geq 4$\\
         \hline

         $m_{ij}$ & $2$ & $3$ & $4$ & $6$ & $\infty$ \\

         \hline

       \end{tabular}

\end{figure}
We also use $r_i$ $(i\in I^{\text{re}})$ to denote the automorphism of $P^\vee$ given by $$r_i(h)=h-\alpha_i(h)h_i \ \ \text{for all}  \ h\in P^\vee.$$
Let $\mathbf i=(i_1,i_2,\dots,i_N)$ be a sequence in $I^{\text{re}}$. Note that for any $\lambda\in \h^*$ and $h\in P^\vee$, we have 
\begin{equation*}
\lambda\big(r_{i_1}r_{i_2}\cdots{r_{i_N}}(h)\big)=\big(r_{i_N}r_{i_{N-1}}\cdots{r_{i_1}}(\lambda)\big)(h).
\end{equation*}

\vskip 3mm

Let $I^\infty=(I^{\text{re}}\times \{1\})\cup(I^{\text{im}}\times \Z_{>0})$. For simplicity, we will often write $i$ instead of $(i,1)$ when $i \in I^{\text{re}}$. Let $q$ be an indeterminate, and set for each $i\in I$$$q_i=q^{s_i},\  q_{(i)}=q^{\frac{(\alpha_i,\alpha_i)}{2}}.$$
For $i\in I^{\text{re}}$ and $n\in \Z_{\geq 0}$, we define
$$[n]_i=\frac{q_{i}^n-q_{i}^{-n}}{q_{i}-q_{i}^{-1}},\quad [n]_i!=\prod_{k=1}^n [k]_i,\quad
{\begin{bmatrix} n \\ k \end{bmatrix}}_i=\frac{[n]_i!}{[k]_i![n-k]_i!}.$$

\vskip 3mm

Let $\mathscr E=\Q(q)\left< e_{il} \mid (i,l) \in I^{\infty} \right>$ be the free associative algebra over $\Q(q)$ generated by the symbols $e_{il}$ for $(i,l)\in I^{\infty}$. By setting $\text{deg} e_{il}= l\alpha_{i} $, $\mathscr E$ becomes a $Q_+$-graded algebra. For a homogeneous element $u$ in $\mathscr E$, we denote by $|u|$ the degree of $u$, and for any $A \subseteq Q_{+}$, set
${\mathscr E}_{A}=\{ x\in {\mathscr E} \mid |x| \in A \}$.

\vskip 3mm

We define a {\it twisted} multiplication on $\mathscr E\otimes\mathscr E$ by
$$(x_1\otimes x_2)(y_1\otimes y_2)=q^{(|x_2|,|y_1|)}x_1y_1\otimes x_2y_2$$
for all homogeneous elements $x_1,x_2,y_1,y_2\in\mathscr E$, and equip $\mathscr E$ with a comultiplication $\rho$ defined by
$$\rho(e_{il})=\sum_{m+n=l}q_{(i)}^{mn}e_{im}\otimes e_{in} \ \ \text{for}  \ (i,l)\in I^{\infty}.$$
Here, we set $e_{i0}=1$, and $e_{il}=0$ for $l<0$.

\vskip 3mm

\begin{proposition}\cite{Bozec2014b,Bozec2014c} \
{\rm For a family $\nu=(\nu_{il})_{(i,l)\in I^{\infty}}$ of non-zero elements in $\Q(q)$, there exists a symmetric bilinear form $\{ \ , \ \} :\mathscr E\times\mathscr E\rightarrow \Q(q)$ such that
\begin{itemize}
\item[(a)] $\{x, y\} =0$ if $|x| \neq |y|$,

\item[(b)] $\{1,1\} = 1$,

\item[(c)] $\{e_{il}, e_{il}\} = \nu_{il}$ for all $(i,l) \in
I^{\infty}$,

\item[(d)] $\{x, yz\} = \{\rho(x), y \otimes z\}$  for all $x,y,z
\in {\mathscr E}$.
\end{itemize}
Here, $\{x_1\otimes x_2,y_1\otimes y_2\}=\{x_1,y_1\}\{x_2,y_2\}$ for any $x_1,x_2,y_1,y_2\in {\mathscr E}$.}
\end{proposition}

\vskip 3mm
Let $\mathcal C_n$ be the set of compositions $\mathbf c$ of $n$, and $e_{i,\mathbf c}=e_{ic_1}\cdots e_{ic_m}$ for each $i\in I^{\text{im}}$, $\mathbf c=(c_1,\cdots,c_m)\in \mathcal C_n$. It is clear that $\{e_{i,\mathbf c}\mid \mathbf c\in \mathcal C_n\}$ form a basis of $\mathscr E_{l\alpha_i}$.

\vskip 3mm

Assume that $i\in I^{\text{re}}$, $j\in I$ and $i\neq j$. Let $m\in \Z_{>0}$, $n\in \Z_{\geq 0}$ with $m>-a_{ij}n$, then for any $\mathbf c\in\mathcal C_n$, the following element of $\mathscr E$ belongs to the radical $\mathcal J$ of the form $\{ \ , \ \}$
\begin{equation}\label{qsr}
\mathcal F_{i,j,n,m,\mathbf c,\pm 1}=\sum_{r+s=m}(-1)^rq_i^{\pm r(-a_{ij}n-m+1)}e_i^{(r)}e_{j,\mathbf c}e_i^{(s)}
\end{equation}
Here, if $j \in I^{\text{re}}$, we set $e_{j,\mathbf c}=e_j^{(n)}$ as the divided power of $e_j$. 

Moreover, if $(i,k),(j,l)\in I^{\infty}$ such that $a_{ij}=0$, one can show that the element $e_{ik}e_{jl}-e_{jl}e_{ik}$ belongs to $\mathcal J$.

\vskip 3mm

From now on, we assume that
\begin{equation} \label{eq:assumption}
\nu_{il} \in 1+q^{-1}\Z_{\geq0}[[q^{-1}]]\ \ \text{for all} \ (i,l)\in I^{\infty}.
\end{equation}
Under this assumption, the bilinear form $\{ \ , \ \}$ is non-degenerate on $\mathscr E(i)=\bigoplus _{l\geq 1} {\mathscr E}_{l \alpha_i}$ for $i\in I^{\rm{im}} \backslash I^{\rm{iso}}$. Moreover, its radical is generated by a simpler set consisting of
$$ \sum_{r+s=1-la_{ij}}(-1)^re_i^{(r)}e_{jl}e_i^{(s)} \ \ \text{for} \ i\in
I^{\text{re}},(j,l)\in I^{\infty} \ \text {and} \ i \neq (j,l),$$
and $e_{ik}e_{jl}-e_{jl}e_{ik}$ for all $(i,k),(j,l)\in I^{\infty}$ with $a_{ij}=0$ (cf. \cite[Proposition 14]{Bozec2014b}).

\vskip 3mm

\begin{definition}Given a Borcherds-Cartan datum $(A, P, P^{\vee}, \Pi, \Pi^{\vee})$, the {\it quantum Borcherds-Bozec algebra} $U_q(\g)$ is the associative algebra over $\Q(q)$ with $\mathbf 1$
generated by the elements $q^h$ $(h\in P^{\vee})$ and $e_{il},
f_{il}$ $((i,l) \in I^{\infty})$, subjecting to
$$
\begin{aligned}
& q^0=\mathbf 1,\quad q^hq^{h'}=q^{h+h'} \ \ \text{for} \ h,h' \in P^{\vee} \\
& q^h e_{jl}q^{-h} = q^{l\alpha_j(h)} e_{jl}, \ \ q^h f_{jl}q^{-h} = q^{-l\alpha_j(h)} f_{jl}\ \ \text{for} \ h \in P^{\vee}, (j,l)\in I^{\infty}, \\
& \sum_{r+s=1-la_{ij}}(-1)^r
{e_i}^{(r)}e_{jl}e_i^{(s)}=0 \ \ \text{for} \ i\in
I^{\text{re}},(j,l)\in I^{\infty} \ \text {and} \ i \neq (j,l), \\
& \sum_{r+s=1-la_{ij}}(-1)^r
{f_i}^{(r)}f_{jl}f_i^{(s)}=0 \ \ \text{for} \ i\in
I^{\text{re}},(j,l)\in I^{\infty} \ \text {and} \ i \neq (j,l), \\
& e_{ik}e_{jl}-e_{jl}e_{ik} = f_{ik}f_{jl}-f_{jl}f_{ik} =0 \ \ \text{for} \ a_{ij}=0,\\
& e_{ik}f_{jl}=f_{jl}e_{ik} \ \ \text{for} \ (i,k),(j,l)\in I^{\infty} \ \text{and}\ i\neq j,\\
& \sum_{\substack{m+n=k \\ n+s=l}}q_{(i)}^{n(m-s)}\nu_{in}K_i^{-n}f_{is}e_{im}=\sum_{\substack{m+n=k \\ n+s=l}}q_{(i)}^{n(s-m)}\nu_{in}K_i^{n} e_{im}f_{is}\ \ \text{for} \ (i,k),(i,l)\in I^{\infty}.
\end{aligned}
$$
Here, $K_i=q_i^{h_i}$ for all $i\in I$. We extend the grading by setting $|q^h|=0$ and $|f_{il}|= -l \alpha_{i}$. For each $\beta=\sum n_i\alpha_i\in Q$, we set $K_\beta=\prod_iK_i^{n_i}$.
\end{definition}

\vskip 3mm

Let $U^+_q(\g)$ (resp. $U^-_q(\g)$) be the subalgebra of $U_q(\g)$ generated by $e_{il}$ (resp. $f_{il}$) for $(i,l)\in I^{\infty}$,
and $U^{0}_q(\g)$ the subalgebra of $U_q(\g)$ generated by $q^h$ for $h\in P^{\vee}$.
Then the quantum Borcherds-Boec algebra $U_q(\g)$ has the {\it triangular decomposition}
$$U_q{(\g)}\cong U^-_q(\g)\otimes U^0_q(\g) \otimes U^+_q(\g).$$
We shall denote by $U$ (resp. $U^+$ and $U^-$) for $U_q(\g)$ (resp. $U^+_q(\g)$ and $U^-_q(\g)$) for simplicity.

\vskip 3mm
The algebra $U$ is endowed with a comultiplication
$\Delta: U \rightarrow U \otimes U$
given by
\begin{equation} \label{eq:comult}
\begin{aligned}
& \Delta(q^h) = q^h \otimes q^h, \\
& \Delta(e_{il}) = \sum_{m+n=l} q_{(i)}^{mn}e_{im}K_{i}^{n}\otimes e_{in}, \\
& \Delta(f_{il}) = \sum_{m+n=l} q_{(i)}^{-mn}f_{im}\otimes K_{i}^{-m}f_{in}.
\end{aligned}
\end{equation}

\vskip 3mm

We shall give the relations between $\rho:U^+ \rightarrow U^+\otimes U^+$ and $\Delta$. Let $x\in U^+$ be a homogeneous element such that $\rho(x)=\sum x_1\otimes x_2$, then we have
\begin{equation}\label{b}
\Delta(x)=\sum x_1K_{|x_2|}\otimes x_2.
\end{equation}
 If $y\in U^-$ is a homogeneous element with $\rho(\omega (y))=\sum y_1\otimes y_2$, where $\omega$  is the involution of $U$ such that
$\omega(q^h)=q^{-h},\ \omega(e_{il})=f_{il},\ \omega(f_{il})=e_{il}$ for $h \in P^{\vee}$ and $(i,l)\in I^{\infty}$. We have
\begin{equation}\label{bb}
\Delta(y)=\sum q^{-(|y_1|,|y_2|)}\omega (y_2)\otimes K_{-|y_2|}\omega (y_1).
\end{equation}



\vskip 3mm

\begin{proposition}\cite{Bozec2014b,Bozec2014c}\label{prim}
{\rm For any $i\in I^{\text {im}}$ and $l\geq 1$, there exist unique elements $\mathtt a_{il}\in  U^+_{l \alpha_{i}}$ and $\mathtt b_{il}=\omega (\mathtt a_{il})$ such that
\begin{itemize}\label{bozec}
\item[(1)] $\Q (q) \left<e_{il} \mid l\geq 1\right>=\Q (q) \left<\mathtt a_{il} \mid l\geq 1\right>$ and $\Q (q) \left<f_{il} \mid l\geq 1\right>=\Q (q) \left<\mathtt b_{il} \mid l\geq 1\right>$,
\item[(2)] $(\mathtt a_{il},z)_L=0$ for all $z\in \Q (q) \left<e_{i1} ,\cdots,e_{il-1}\right>$,\\
$(\mathtt b_{il},z)_L=0$ for all $z\in \Q (q) \left<f_{i1} ,\cdots,f_{il-1}\right>$,
\item[(3)] $\mathtt a_{il}-e_{il}\in \Q(q) \left<e_{ik} \mid k<l \right>$ and $\mathtt b_{il}-f_{il}\in \Q (q) \left<f_{ik} \mid k<l \right>$,
\item[(4)] $\overline{\mathtt a}_{il}=\mathtt a_{il},\ \overline{\mathtt b}_{il}=\mathtt b_{il}$,
\item[(5)] $\rho(\mathtt a_{il})=\mathtt a_{il}\otimes 1+1\otimes \mathtt a_{il}, \ \rho(\mathtt b_{il})=\mathtt b_{il}\otimes 1+ 1\otimes \mathtt b_{il}$.
\end{itemize}
Here, $^-:U^{\pm}\rightarrow U^{\pm}$ is the $\Q$-algebra homomorphism defined by
$\overline{e}_{il}=e_{il},\ \overline{f}_{il}=f_{il}$ and  $\overline{q}=q^{-1}$.
\
}
\end{proposition}

Set $\tau_{il}=\{\mathtt a_{il},\mathtt a_{il}\}=\{\mathtt b_{il},\mathtt b_{il}\}$. We have the following commutation relations in $U$ derived from the Drinfeld double process,
\begin{equation}\label{news}
\mathtt a_{il}\mathtt b_{jk}-\mathtt b_{jk}\mathtt a_{il}=\delta_{ij}\delta_{lk}\tau_{il}(K_i^{-l}-K_i^{l}).
\end{equation}
The $\mathtt a_{il}$'s and $\mathtt b_{il}$'s are called the {\it primitive generators} of $U_q(\g)$.

\vskip 3mm

Let $\mathcal C_{l}$ (resp. $\mathcal P_{l}$) be the set of compositions (resp. partitions) of $l$. For $i\in I^{\text{im}}$, we define
\begin{equation*}
\mathcal D_{i,l}=\begin{cases} \mathcal C_{l}& \text{if} \ i\in I^{\text{im}}\backslash I^{\text{iso}}, \\ \mathcal P_{l}& \text{if} \ i\in I^{\text{iso}}.\end{cases}
\end{equation*}
and $\mathcal D_i=\bigsqcup_{l\geq 0}\mathcal D_{i,l}$. Let $\mathbf c =(c_1,\cdots,c_t)\in \mathcal D_{i,l}$,  we set
  $$\mathtt a_{i,\mathbf c}=\mathtt a_{ic_1}\cdots\mathtt a_{ic_t}, \  \mathtt b_{i,\mathbf c}=\mathtt b_{ic_1}\cdots\mathtt b_{ic_t}  \ \text{and} \ \tau_{i,\mathbf c}=\tau_{ic_1}\cdots\tau_{ic_t} .$$
 Note that $\{ \mathtt a_{i,\mathbf c} \mid \mathbf c \in \mathcal D_{i,l}\}$ forms a basis of  $U^+_{l \alpha_{i}}$. For each $i\in I^{\text{re}}$, we set $\mathtt a_{i1}=e_{i1}$, $\mathtt b_{i1}=f_{i1}$, and
write $\mathtt a_i$ (resp. $\mathtt b_i$) instead of $\mathtt a_{i1}$ (resp. $\mathtt b_{i1}$) in this case for simplicity.

\vskip 3mm
\begin{example}{\rm
$\lambda\in\mathcal P_l$ can be written as the form $\lambda=1^{\lambda_1}2^{\lambda_2}\cdots l^{\lambda_l}$, where $\lambda_k$ are non-negative integers such that $\lambda_1+2\lambda_2+\cdots+l\lambda_l=l$. For $i\in I^{\text{iso}}$, we have
$$\mathtt a_{il}=e_{il}-\sum_{\lambda\in \mathcal P_l\backslash (l)} \frac{1}{\prod_{i=1}^{l} \lambda_i!} \mathtt a_{i,\lambda}.$$
}\end{example}
\vskip 3mm

\begin{example}{\rm
 For $i\in I^{\text{iso}}$ and $\mathbf c, \mathbf c' \in \mathcal P_{l}$,  we have
$\{\mathtt a_{i,\mathbf c},\mathtt a_{i,\mathbf c'}\}=0$ when $\mathbf c\neq \mathbf c' $.
 For $i\in I^{\text{im}}\backslash I^{\text{iso}}$ and $\mathbf c, \mathbf c' \in \mathcal C_{l}$, if the partitions obtained by rearranging $\mathbf c$ and $\mathbf c'$  are not equal, we also have 
$\{\mathtt a_{i,\mathbf c},\mathtt a_{i,\mathbf c'}\}=\{\mathtt b_{i,\mathbf c},\mathtt b_{i,\mathbf c'}\}=0$.}
\end{example}

\vskip 3mm

\begin{example}{\rm
Let $\mathbf c =(c_1,c_2,\cdots,c_r)$ be a composition of $l$, we denote by $\widetilde{\mathbf c}$ the reverse of $\mathbf c$, i.e. $\widetilde{\mathbf c}=(c_r,c_{r-1},\cdots,c_1)$. Assume that $\mathbf c,\mathbf c'$ are two compositions of $l$, which determine the same partition, say it is $\mathbf p =(p_1,p_2,\cdots,p_r)$. For $i\in I^{\text{im}} \backslash I^{\text{iso}}$, if $\{\mathtt a_{i,\mathbf c},\mathtt a_{i,\mathbf c'}\}=\gamma \tau_{i,\mathbf p}$ for some $\gamma \in \Q(q)$, then we have
$$\{\mathtt a_{i,\mathbf c},\mathtt a_{i,\widetilde{\mathbf c}'}\}=q_{(i)}^{m} \overline{\gamma} \tau_{i,\mathbf p}.$$

Here, $m=2\sum_{r<s}c_rc_s=2\sum_{r<s}c_r'c_s'$ and $^-$ be the involution of $\Q(q)$ mapping $q$ to $q^{-1}$. In particular, we see that
$\{\mathtt a_{i,\widetilde{\mathbf c}},\mathtt a_{i,\widetilde{\mathbf c}'}\}=\{\mathtt a_{i,\mathbf c},\mathtt a_{i,\mathbf c'}\}$.
}\end{example}

\vskip 3mm

\begin{definition}
For every $(i,l)\in I^{\infty}$, we define the linear maps $\delta_{i,l}, \delta^{i,l}:U^{+}\rightarrow U^{+}$ by
$$\delta_{i,l}(1)=0,\ \delta_{i,l}(\mathtt a_{jk})=\delta_{ij}\delta_{lk}\ \text{and} \ \delta_{i,l}(xy)=q^{l(|y|,\alpha_i)}\delta_{i,l}(x)y+x\delta_{i,l}(y),$$
$$\delta^{i,l}(1)=0,\ \delta^{i,l}(\mathtt a_{jk})=\delta_{ij}\delta_{lk}\ \text{and} \ \delta^{i,l}(xy)=\delta^{i,l}(x)y+q^{l(|x|,\alpha_i)}x\delta^{i,l}(y),$$
for any homogeneous elements $x,y$ in $U^+$.
\end{definition}

\vskip 3mm

Let $(i,l) \in I^{\infty}$ and $z\in U^+$, one deduce from the Drinfeld double process that (cf. \cite[Proposition 3.10]{Bozec2014c})
\begin{equation}\label{c}
 \left[\mathtt a_{il},\omega(z)\right]=\tau_{il}\{\omega(\delta_{i,l}(z))K_i^{-l}-K_i^{l}\omega(\delta^{i,l}(z))\}.
\end{equation}
Applying the involution $\omega$ to both sides of (\ref{c}), we obtain the following commutation relation
\begin{equation}
 \left[\mathtt b_{il},z \right]=\tau_{il}\{ \delta_{i,l}(z)K_i^{l}-K_i^{-l}\delta^{i,l}(z)\}.
 \end{equation}

\vskip 8mm
\section{Lusztig symmetries and the braid group actions}

\vskip 2mm

Let $U_q(\g)$ be a quantum Borcherds-Bozec algebra, we denote by $\mathcal O$ the category of $U_q(\g)$-module $M$ satisfying the following two conditions \begin{itemize}
\item[(i)] $M$ has a weight space decomposition
$$M=\bigoplus_{\mu\in P}M_\mu, \ \text{where}  \ M_\mu=\{z\in M\mid q^hz=q^{\mu(h)}z\},$$
\item[(ii)] $\mathtt a_i$ and $\mathtt b_i$ act locally nilpotent on $M$ for all $i\in I^{\text{re}}$.
\end{itemize}

\vskip 3mm

For $i\in I^{\text{re}}$, we set $B_i=\mathtt b_i \big/ \tau_i(q_i^{-1}-q_i)$, then the commutation relation between $\mathtt a_i$ and $\mathtt B_i$ is given by
$$\mathtt a_i\mathtt B_i-\mathtt B_i\mathtt a_i={(K_i-K_i^{-1})}/{(q_i-q_i^{-1})}.$$
Since every element $M$ in $\mathcal O$ is a weight module by (i), so it has a direct sum decomposition as follow
$$M=\bigoplus_{n\in \Z}M^n, \ \text{where}  \ M^n=\{z\in M\mid K_iz=q_i^nz\}.$$

We see that $\mathtt a_i(M^n)\subseteq M^{n+2}$, $\mathtt B_i(M^n)\subseteq M^{n-2}$ for all $n$, and $\mathtt a_i\mathtt B_i-\mathtt B_i\mathtt a_i$ acts  on $M^n$ is the multiplication by $[n]_i$ for all $n$, these properties yield $M$ becomes an object of the category $\mathcal C'_i$ defined in \cite[5.1.1]{Lusztig}. Hence the standard argument in \cite[Chapter 5]{Lusztig} can be applied to our case. We define the $\Q(q)$-linear maps $L'_{i,e},L''_{i,e}:M\rightarrow M$ $(e=\pm 1)$ for $i\in I^\text{re}$ and $M\in \mathcal O$ by
$$L'_{i,e}(z)=\sum_{a,b,c;a-b+c=n}(-1)^bq_i^{e(-ac+b)}\mathtt B_i^{(a)}\mathtt a_i^{(b)}\mathtt B_i^{(c)}z,$$
$$L''_{i,e}(z)=\sum_{a,b,c;-a+b-c=n}(-1)^bq_i^{e(-ac+b)}\mathtt a_i^{(a)}\mathtt B_i^{(b)}\mathtt a_i^{(c)}z,$$
where $z\in M^n$ and $e=\pm 1$. $L'_{i,e},L''_{i,e}$ are called {\it symmetries} on $M$. According to \cite[5.2.3]{Lusztig}, we have
$L'_{i,e}L''_{i,-e}=L''_{i,-e}L'_{i,e}={\text{id}}:M\rightarrow M$.

\vskip 3mm

Given $i\in I^{\text{re}}$, $j\in I$ and $i\neq j$. Let $m,n \in \Z_{\geq 0}$. Along the notations in (\ref{qsr}), we set
\begin{equation}\label{F}
\mathcal F_{i,j,n,m,e}=\sum_{r+s=m}(-1)^rq_i^{er(-a_{ij} n-m+1)}\mathtt a_i^{(r)}\mathtt a_{jn}\mathtt a_i^{(s)}\in U^+.
\end{equation}
Here $a_{jn}=1$ if $n=0$. Otherwise, $\mathtt a_{jn}$ are the generators of $U^+$ corresponding to $(j,n)\in I^{\infty}$ when $j\in I^{\text{im}}$, and $\mathtt a_{jn}=a_j^{(n)}$ when $j\in I^{\text{re}}$.

We shall denote $\mathcal F_{i,j,n,m,e}$ by $\mathcal F_{n,m,e}$  for simplicity if there is no risk of confusion, and set $\beta=-a_{ij}$, $\beta'=-a_{ji}$.

\vskip 3mm

We shall give several equations about $\mathcal F_{n,m,e}$ that can be proved by similar inductive processes in \cite[Chapter 7]{Lusztig}.

\vskip 3mm

\begin{lemma}\label{eqq}\
{\rm
\begin{itemize}
\item [(i)] $-q_i^{e(\beta n-2m)}\mathtt a_i\mathcal F_{n,m,e}+\mathcal F_{n,m,e}\mathtt a_i={[m+1]}_i\mathcal F_{n,m+1,e}$,\\
\item [(ii)] $\mathtt a_i^{(p)}\mathcal F_{n,m,e}=\sum_{p'=0}^p(-1)^{p'}q_i^{e(2pm-\beta pn+pp'-p')}{\left[ \substack{ m+p'\\p'}\right]}_i\mathcal F_{n,m+p',e}\mathtt a_i^{(p-p')}$,\\
\item [(iii)] $-\mathtt B_i\mathcal F_{n,m,e}+\mathcal F_{n,m,e}\mathtt B_i={[\beta n-m+1]}_iK_{-ei}\mathcal F_{n,m-1,e}$,\\
\item [(iv)] $\mathtt B_i^{(p)}\mathcal F_{n,m,e}=\sum_{p'=0}^p(-1)^{p'}q_i^{-e(pp'-p')}{\left[ \substack{ \beta n-m+p'\\p'}\right]}_iK_{-ep'i}\mathcal F_{n,m-p',e}\mathtt B_i^{(p-p')}$.
\end{itemize}
}\end{lemma}

\vskip 3mm

\begin{lemma}\ {\rm
\begin{itemize}
\item [(i)] If $j\in I^{\text{re}}$ and $m=1+\beta n$, then
$$\mathtt B_j\mathcal F_{n,m,e}-\mathcal F_{n,m,e}\mathtt B_j=K_{-j}\frac{q_j^{n-1}}{q_j-q_j^{-1}}\mathcal F_{n-1,m,1}-K_{j}\frac{q_j^{1-n}}{q_j-q_j^{-1}}\mathcal F_{n-1,m,-1}.$$
\item [(ii)] If $j\in I^{\text{im}}$, $(j,l)\in I^\infty$  and $m=1+\beta n$, then
$$\mathtt b_{jl}\mathcal F_{n,m,e}-\mathcal F_{n,m,e}\mathtt b_{jl}=0.$$
\end{itemize}
\begin{proof}
In the case of $j\in I^{\text{re}}$, the proof is the same as \cite[Lemma 7.1.4]{Lusztig}. We now prove (ii). Note that
$$\begin{aligned}
& \mathtt b_{jl}\mathcal F_{n,m,e}-\mathcal F_{n,m,e}\mathtt b_{jl}= \sum_{r+s=m}(-1)^rq_i^{er(\beta n-m+1)}\mathtt a_i^{(r)}\big(\mathtt b_{jl}\mathtt a_{jn}-\mathtt a_{jn}\mathtt b_{jl}\big)\mathtt a_i^{(s)}\\
& \phantom{\mathtt b_{jl}\mathcal X_{n,m,e}}=\delta_{ln}\tau_{jn} \sum_{r+s=m}(-1)^rq_i^{er(\beta n-m+1)}\mathtt a_i^{(r)}(K_j^n-K_j^{-n})\mathtt a_i^{(s)}\\
& \phantom{\mathtt b_{jl}\mathcal X_{n,m,e}}=\delta_{ln}\tau_{jn} \sum_{r+s=m}(-1)^rq_i^{er(\beta n-m+1)}(q_i^{nr\beta}K_j^n-q_i^{-nr\beta}K_j^{-n})\mathtt a_i^{(r)}\mathtt a_i^{(s)}.\\
\end{aligned}$$
Since $m=1+\beta n$, we see that the right hand side of the above equality equals to
$$\delta_{ln}\tau_{jn} \sum_{r+s=m}(-1)^r(q_i^{r(m-1)}K_j^n-q_i^{r(1-m)}K_j^{-n})\frac{1}{[r]_i![m-r]_i!}\mathtt a_i^{m},$$
and our assertion follows by the identity
$\sum_{r=0}^{m}(-1)^rq_i^{\pm r(1-m)}{\begin{bmatrix}m\\ r\end{bmatrix}}_i=0$
for all $m>0$.
\end{proof}
}\end{lemma}

\vskip 3mm

 The following lemma is an analogue of \cite[Lemma 3.5.4]{Lusztig}. One just need  to note the fact that if $i\in I^{\text{re}}$, $(j,l)\in I^\infty$ and $i\neq j$, then
$$\mathtt a_i^m\mathtt a_{jl}\in\sum_{k=0}^{-la_{ij}}\Q(q)\mathtt a_i^k\mathtt a_{jl}\mathtt a_i^{m-k}$$
for all $m\geq 1-la_{ij}$.

\vskip 3mm

\begin{lemma}\label{ann}{\rm
Assume that $I^{\text{re}}\neq \emptyset$. Let $u$ be an element of $U$ such that $u$ annihilates all $M\in \mathcal O$. Then $u=0$.
}\end{lemma}

\vskip 3mm

Given an $i\in I^{\text{re}}$, we define the symmetries $L'_{i,e}:U\rightarrow U$ $(e=\pm 1)$ on the generators of $U$ as follows
$$\begin{aligned}
& L'_{i,e}(\mathtt a_i)=-K_{ei}\mathtt B_i,\ L'_{i,e}(\mathtt B_i)=-\mathtt a_i K_{-ei};\\
& L'_{i,e}(\mathtt a_{jl})=\sum_{r+s=-la_{ij}}(-1)^rq_i^{er}\mathtt a_i^{(r)}\mathtt a_{jl}\mathtt a_i^{(s)} \ \ \text{for}\ (j,l)\in I^\infty \ \text {and} \  i\neq (j,l);\\
& L'_{i,e}(\mathtt b_{jl})=\sum_{r+s=-la_{ij}}(-1)^rq_i^{-er}\mathtt B_i^{(s)}\mathtt b_{jl}\mathtt B_i^{(r)} \ \ \text{for}\ (j,l)\in I^\infty \ \text {and} \  i\neq (j,l);\\
& L'_{i,e}(q^h)=q^{r_i(h)}=q^{h-\alpha_i(h)h_i },
\end{aligned}$$
and define the symmetries $L''_{i,-e}$ $(e=\pm 1)$ by
$$\begin{aligned}
& L''_{i,-e}(\mathtt a_i)=-\mathtt B_iK_{-ei},\ L''_{i,-e}(\mathtt B_i)=-K_{ei}\mathtt a_i ;\\
& L''_{i,-e}(\mathtt a_{jl})=\sum_{r+s=-la_{ij}}(-1)^rq_i^{er}\mathtt a_i^{(s)}\mathtt a_{jl}\mathtt a_i^{(r)} \ \ \text{for}\ (j,l)\in I^\infty \ \text {and} \  i\neq (j,l);\\
& L''_{i,-e}(\mathtt b_{jl})=\sum_{r+s=-la_{ij}}(-1)^rq_i^{-er}\mathtt B_i^{(r)}\mathtt b_{jl}\mathtt B_i^{(s)} \ \ \text{for}\ (j,l)\in I^\infty \ \text {and} \  i\neq (j,l);\\
& L''_{i,-e}(q^h)=q^{h-\alpha_i(h)h_i }.
\end{aligned}$$

\vskip 3mm

Let $j\in I$ with $j\neq i$, and let $m,n \in \Z_{\geq 0}$. We set $\mathcal F_{i,j,n,m,e}=\mathcal F_{n,m,e}$ as in (\ref{F}) and 
\begin{equation}\label{FG}
 \begin{aligned}
& \mathcal F'_{i,j,n,m,e}=\mathcal F'_{n,m,e}=\sum_{r+s=m}(-1)^rq_i^{er(\beta n-m+1)}\mathtt a_i^{(s)}\mathtt a_{jn}\mathtt a_i^{(r)};\\
& \mathcal G_{i,j,n,m,e}=\mathcal G_{n,m,e}=\sum_{r+s=m}(-1)^rq_i^{-er(\beta n-m+1)}\mathtt B_i^{(s)}\mathtt b_{jn}\mathtt B_i^{(r)};\\
& \mathcal G'_{i,j,n,m,e}=\mathcal G'_{n,m,e}=\sum_{r+s=m}(-1)^rq_i^{-er(\beta n-m+1)}\mathtt B_i^{(r)}\mathtt b_{jn}\mathtt B_i^{(s)}.
\end{aligned}
\end{equation}
Here, $\mathtt b_{jn}$ has the same interpretation as $\mathtt a_{jn}$ and $\beta=-a_{ij}$.

\vskip 3mm

 Using Lemma \ref{eqq} and Lemma \ref{ann}, we can verify the following statements step by step according to \cite[Chapter 37]{Lusztig}, we leave it to readers. Precisely, thanks to Lemma \ref{eqq}, one can treat the generator $\mathtt a_{jn}$  as the $\mathtt a_j^{(n)}$ in form when we work with (i) in the following proposition, which is a counterpart of Lemma 37.2.2 in {\it loc.cit.}. Then the proof of (ii) is entirely similar to  Lemma 37.2.3 in \cite{Lusztig} with the help of (i) and Lemma \ref{ann}.

\vskip 3mm

\begin{proposition}\label{Main}\
{\rm
\begin{itemize}
\item [(i)] Let $M\in \mathcal O$ and $z\in M$, then we have
$$\begin{aligned}
& L''_{i,-e}(\mathcal F_{n,\beta n,e}z)=\mathtt a_{jn}L''_{i,-e}(z);\ \ L'_{i,e}(\mathcal F'_{n,\beta n,e}z)=\mathtt a_{jn}L'_{i,e}(z);\\
& L''_{i,-e}(\mathcal G_{n,\beta n,e}z)=\mathtt b_{jn}L''_{i,-e}(z);\ \ L'_{i,e}(\mathcal G'_{n,\beta n,e}z)=\mathtt b_{jn}L'_{i,e}(z).
\end{aligned}$$
\item [(ii)] $L'_{i,e}$ and $L''_{i,-e}$ are automorphisms of $U$ sending the $\alpha$-root space onto the $r_i(\alpha)$-root space, they are the inverse of each other. For any $M\in \mathcal O$, let $z\in M$ and $u\in U$, we have
$$L'_{i,e}(uz)=L'_{i,e}(u)L'_{i,e}(z)\ \text{and} \ L''_{i,-e}(uz)=L''_{i,-e}(u)L''_{i,-e}(z).$$
Furthermore, the operators $L'_{i,e}$ and $L''_{i,-e}$ of $U$ are uniquely determined by these properties.
\item [(iii)] For any $m\in \Z$, we have
$$L'_{i,e}(\mathcal F'_{n,m,e})=\mathcal F_{n,\beta n-m,e}\ \text{and} \ L''_{i,-e}(\mathcal F_{n,m,e})=\mathcal F'_{n,\beta n-m,e}.$$
\end{itemize}
}\end{proposition}

\vskip 3mm

It is more convenient to twist the involution $\omega$. We shall denote by $\varpi$ the automorphism of $U$ such that
$$\varpi(\mathtt a_i)=\mathtt B_i,\ \varpi(\mathtt B_i)=\mathtt a_i,\ \varpi(q^h)=q^{-h}, \ \varpi(\mathtt a_{jl})=\mathtt b_{jl}\ \text{and} \ \varpi(\mathtt b_{jl})=\mathtt a_{jl}\ \text{for}\ (j,l)\neq i.$$
Let $*$ be the anti-automorphism of $U$ given by
$$*(q^h)=q^{-h},\ *(\mathtt a_{jl})=\mathtt a_{jl}\ \text{and}\ *(\mathtt b_{jl})=\mathtt b_{jl}\ \ \text{for}\ h \in P^{\vee},\ (i,l)\in I^{\infty}.$$

Checking for the generators of $U$, we obtain $\varpi L'_{i,e}=L''_{i,e}\varpi: U\rightarrow U$ and $*L'_{i,e}=L''_{i,-e}*:U\rightarrow U$.
Moreover, if $u\in U$ such that $K_iuK_i^{-1}=q_i^nu$, we have \begin{equation}\label{L'L''}L''_{i,e}(u)=(-1)^nq_i^{en}L'_{i,e}(u).\end{equation}

\vskip 3mm

\begin{lemma}\cite[Lemma 39.4.1]{Lusztig} {\rm Assume that $I=I^{\text{re}}$ consists of exact two elements $i\neq j$ such that $m_{ij}=2,3,4$ or $6$, and assume $M\in \mathcal O$ in this case,  we have
\begin{equation}\label{L0}
L''_{i,-e}L''_{j,-e}L''_{i,-e}\cdots=L''_{j,-e}L''_{i,-e}L''_{j,-e}\cdots:M\rightarrow M
\end{equation}
with both sides have $m_{ij}$ terms. 
}\end{lemma}

\vskip 1mm
We shall generalize it to an arbitrary Borcherds-Cartan datum by the meothod in \cite[Lemma 39.4.3]{Lusztig}.

\vskip 3mm

\begin{theorem}\label{Ba}{\rm
Given a  Borcherds-Cartan datum $(A, P, P^{\vee} , \Pi, \Pi^{\vee})$ with $I^{\text{re}}\neq \emptyset$, let $U$ be the associated quantum Borcherds-Bozec algebra. For any $i\neq j$ such that $m_{ij}<\infty$, we have the following equalities of automorphisms of $U$ and any $U$-module $M$ in $\mathcal O$
\begin{equation}\label{L1}
 L''_{i,-e}L''_{j,-e}L''_{i,-e}\cdots=L''_{j,-e}L''_{i,-e}L''_{j,-e}\cdots ,
\end{equation}
\begin{equation}\label{L2}
 L'_{i,e}L'_{j,e}L'_{i,e}\cdots=L'_{j,e}L'_{i,e}L'_{j,e}\cdots,
\end{equation}
where all the products have $m_{ij}$ factors.
\begin{proof}
By the preceding argument,  equality (\ref{L1}) for $M\in \mathcal O$ can be shown by considering $M$ as a $U'$-module, where $U'$ is the subalgebra of $U$ generated by $\mathtt a_i,\mathtt b_i, \mathtt a_j, \mathtt b_j$ and $U^0$.

Let $u\in U$,  we set $u_1=(L''_{i,-e}L''_{j,-e}L''_{i,-e}\cdots)(u)$ and $u_2=(L''_{j,-e}L''_{i,-e}L''_{j,-e}\cdots)(u)$.For any $M\in \mathcal O$ and $z\in M$, we have by using Proposition \ref{Main} (ii) that
$$\begin{aligned}
& u_1((L''_{j,-e}L''_{i,-e}L''_{j,-e}\cdots)z)=u_1((L''_{i,-e}L''_{j,-e}L''_{i,-e}\cdots)z)\\
& =(L''_{i,-e}L''_{j,-e}L''_{i,-e}\cdots)(uz)=(L''_{j,-e}L''_{i,-e}L''_{j,-e}\cdots)(uz)\\
& =u_2((L''_{j,-e}L''_{i,-e}L''_{j,-e}\cdots)z).
\end{aligned}$$
It follows that $u_1-u_2$ annihilates all $M\in \mathcal O$. Hence we get $u_1=u_2$ by Lemma \ref{ann}, this proves (\ref{L1}) for $U$. Finally, the equality for $L'$ follows by taking inverse.
\end{proof}
}\end{theorem}

\vskip 3mm

The {\it braid group} associated to a Borcherds-Cartan datum is the group generated by $r_i$ $(i\in I^{\text{re}})$ with defining relations
$$r_ir_jr_i\cdots=r_jr_ir_j\cdots  \quad \text{for}\ i\neq j \ \text{and} \ m_{ij}<\infty,$$
where both sides have $m_{ij}$ factors. Theorem \ref{Ba} yields the braid group actions on $U$ and $U$-module $M\in \mathcal O$.

\vskip 3mm

If $r_{i_1}r_{i_2}\cdots r_{i_N}$ and $r_{i'_1}r_{i'_2}\cdots r_{i'_N}$ are two reduced expressions of $r\in W$, then we have the equality $r_{i_1}r_{i_2}\cdots r_{i_N}=r_{i'_1}r_{i'_2}\cdots r_{i'_N}$ in the braid group.  Hence the following definition is valid
$$L''_{r,e}=L''_{i_1,e}L''_{i_2,e}\cdots L''_{i_N,e},\ L'_{r,e}=L'_{i_1,e}L'_{i_2,e}\cdots L'_{i_N,e},$$
where $r=r_{i_1}r_{i_2}\cdots r_{i_N}$ is a reduced expression of $r\in W$. From the definition, we have
$$L''_{rr',e}=L''_{r,e}L''_{r',e},\ L'_{rr',e}=L'_{r,e}L'_{r',e}$$
if $r,r'\in W$ such that $l(rr')=l(r)+l(r')$. The almost same calculate in \cite[Chapter 40]{Lusztig} gives the following lemma.

\begin{lemma}{\rm
Let $r\in W$ and let $i\in I^{\text{re}}$ with $l(rr_i)=l(r)+1$, we have $L''_{r,e}(\mathtt a_i)\in U^+$ and $L'_{r,e}(\mathtt a_i)\in U^+$ for $e=\pm 1$.
}\end{lemma}

\vskip 3mm

\begin{corollary}{\rm
Let $r_{i_1}r_{i_2}\cdots r_{i_N}$ be a reduced expression for some $r\in W$. Then
\begin{itemize}
\item [(i)] $L''_{i_1,e}L''_{i_2,e}\cdots L''_{i_{N-1},e}(\mathtt a_{i_N})\in U^+$ and $L'_{i_1,e}L'_{i_2,e}\cdots L'_{i_{N-1},e}(\mathtt a_{i_N})\in U^+.$
\item [(ii)] For any $j\in I^{\text{im}}$, we have
$L''_{r,e}(\mathtt a_{jl})\in U^+$ and  $L'_{r,e}(\mathtt a_{jl})\in U^+$.
\end{itemize}
\begin{proof}
The first assertion follows from the previous lemma directly. We prove (ii) by induction on $l(r)$. Let $r'=r_{i_1}r_{i_2}\cdots r_{i_{N-1}}$ and $i=i_N$, we have $r=r'r_{i}$. Since $l(r')=l(r)-1$, we have $L''_{r,e}=L''_{r',e}L''_{i,e}$ and $L''_{r',e}(\mathtt a_{jl})\in U^+$ by the induction  hypothesis. Applying $L''_{r,e}$ to $\mathtt a_{jl}$, we have
$$\begin{aligned}
& L''_{r,e}(\mathtt a_{jl})=L''_{r',e}L''_{i,e}(\mathtt a_{jl})=L''_{r',e}(\sum_{r+s=-la_{ij}}(-1)^rq_i^{-er}\mathtt a_i^{(s)}\mathtt a_{jl}\mathtt a_i^{(r)})\\
& \phantom{L''_{r,e}(\mathtt a_{jl})}=\sum_{r+s=-la_{ij}}(-1)^rq_i^{-er}L''_{r',e}(\mathtt a_i^{(s)})L''_{r',e}(\mathtt a_{jl})L''_{r',e}(\mathtt a_i^{(r)}).
\end{aligned}$$
Hence by (i), we deduce that $L''_{r,e}(\mathtt a_{jl})\in U^+$. Similarly, we have $L'_{r,e}(\mathtt a_{jl})\in U^+$.
\end{proof}
}\end{corollary}
\vskip 8mm

\section{Link with the bilinear form on $U$}

\vskip 2mm

Fix $i\in I^{\text{re}}$. For any $(j,l)\in I^\infty$ with $i\neq (j,l)$ and any $m\in \Z$, we set
\begin{itemize}
\item [(i)] $f(i,(j,l),m)=\sum_{k=0}^m(-1)^kq_i^{k(la_{ij}+m-1)}\mathtt a_i^{(k)}\mathtt a_{jl}\mathtt a_i^{(m-k)};$\\
\item [(ii)] $f'(i,(j,l),m)=\sum_{k=0}^m(-1)^kq_i^{k(la_{ij}+m-1)}\mathtt a_i^{(m-k)}\mathtt a_{jl}\mathtt a_i^{(k)};$\\
\item [(iii)] $g(i,(j,l),m)=\sum_{k=0}^m(-1)^kq_i^{k(-la_{ij}-m+1)}\mathtt B_i^{(m-k)}\mathtt b_{jl}\mathtt B_i^{(k)};$\\
\item [(iv)] $g'(i,(j,l),m)=\sum_{k=0}^m(-1)^kq_i^{k(-la_{ij}-m+1)}\mathtt B_i^{(k)}\mathtt b_{jl}\mathtt B_i^{(m-k)}.$
\end{itemize}

\vskip 3mm

Note that $f(i,(j,l),m)=\mathcal F_{i,j,l,m,-1}$, $f'(i,(j,l),m)=\mathcal F'_{i,j,l,m,-1}$, $g(i,(j,l),m)=\mathcal G_{i,j,l,m,-1}$ and $g'(i,(j,l),m)=\mathcal G'_{i,j,l,m,-1}$ by using the notations in (\ref{FG}). Sometimes, we will simply write $f_m, f'_m, g_m, g'_m$ for convenience.

\vskip 3mm

 Let $U^+[i]$ (resp. $U^+_*[i]$, $U^-[i]$ and $U^-_*[i]$) be the subalgebra of $U$ generated by the elements $f(i,(j,l),m)$ (resp. $f'(i,(j,l),m)$, $g(i,(j,l),m)$ and $g'(i,(j,l),m)$) for all $(j,l)\neq i$ and $m\in \Z$. We have $*(U^+[i])=U^+_*[i]$ and $*(U^-[i])=U^-_*[i]$.

\vskip 3mm

According to Proposition \ref{Main}(iii) and equality (\ref{L'L''}), for any $m\in \Z$, we have $$L''_{i,1}(f_m)=(-1)^nq_i^nL'_{i,1}(f_m)=f'_{l\beta-m},$$ where $n=2m-l\beta$, and we obtain by using $*$ that $$L'_{i,-1}(f'_m)=(-1)^nq_i^nL''_{i,-1}(f'_m)=f_{l\beta-m}.$$  

Note that
$$\varpi f_m=(-1)^mq_i^{m(la_{ij}+m-1)}g_m,\ \varpi f'_m=(-1)^mq_i^{m(la_{ij}+m-1)}g'_m,$$
so we obtain
$$\begin{aligned}
& L''_{i,1}(g_m)=(-1)^mq_i^{m(l\beta-m+1)}L''_{i,1}(\varpi f_m)\\
&\phantom{L''_{i,1}}=(-1)^mq_i^{m(l\beta-m+1)}\varpi L'_{i,1} (f_m)\\
&\phantom{L''_{i,1}}=(-1)^{m+n}q_i^{m(l\beta-m+1)-n}\varpi f'_{l\beta-m}=g'_{l\beta-m},\\
\end{aligned}$$
and $L'_{i,-1}(g'_m)=g_{l\beta-m}$ by taking inverse. Hence we have $$\varpi U^+[i]=U^-[i], \ \varpi U^+_*[i]=U^-_*[i],$$ and $$L''_{i,1}U^+[i]=U^+_*[i]; \ L''_{i,1}U^-[i]=U^-_*[i];$$
$$L'_{i,-1}U^+_*[i]=U^+[i]; \ L'_{i,-1}U^-_*[i]=U^-[i].$$

\vskip 3mm

\begin{lemma}\label{decom}{\rm We have
\begin{itemize}
\item [(i)]  $U^+=\bigoplus_{t\geq 0}\mathtt a_i^tU^+[i]=\bigoplus_{t\geq 0}U^+[i]\mathtt a_i^t$;
\item [(ii)]  $U^+=\bigoplus_{t\geq 0}\mathtt a_i^t U^+_*[i]=\bigoplus_{t\geq 0} U^+_*[i]\mathtt a_i^t$.
\end{itemize}
\begin{proof}
Note that (ii) follows from (i) by applying $*$. To prove (i), it suffices to show that $U^+[i]=\text{ker} \delta^i$.  By Lemma \ref{eqq}(i), we see that $\mathtt a_i$ and $f_m$ interchange by the following formula
$$-q_i^{2m-l\beta}\mathtt a_i f(i,(j,l),m)+f(i,(j,l),m)\mathtt a_i={[m+1]}_if(i,(j,l),m+1).$$
Note that $U^+$ is generated by $\mathtt a_i$ and $U^+[i]$, we have the following decomposition $$U^+=\sum_{t\geq 0}\mathtt a_i^tU^+[i].$$

Since $L''_{i,1}U^+[i]\subseteq U^+$, one can show by the exposition in  \cite[Lemma 38.1.4]{Lusztig} that $U^+[i]$ is contained in $ \text{ker} \delta^i$. Let $x\in \text{ker}\delta^i$, we write $x$ into the form $x=\sum_{t\geq0}\mathtt a_i^tx_t$, where $x_t\in U^+[i]\subseteq \text{ker} \delta^i$. Note that $U^+=\bigoplus_{t\geq 0}\mathtt a_i^t\text{ker} \delta^i$, we have $x=x_0\in U^+[i]$. This completes the proof.
\end{proof}
}\end{lemma}

\vskip 3mm

Using $\varpi$, we have the decompositions for $U^-$
$$U^-=\bigoplus_{t\geq 0}\mathtt B_i^tU^-_*[i]=\bigoplus_{t\geq 0}U^-_*[i]\mathtt B_i^t=\bigoplus_{t\geq 0}\mathtt B_i^t U^-_*[i]=\bigoplus_{t\geq 0} U^-_*[i]\mathtt B_i^t.$$

\vskip 3mm
We now assume that $x\in U^+$ and $L''_{i,1}(x)\in U^+$, and write $x=\sum_{\alpha\in Q^+}x_\alpha$.  Since $L''_{i,1}$ sends the different root spaces into the different root spaces, we have $L''_{i,1}(x_\alpha)\in U^+$ for each $\alpha\in Q^+$. Note that $L''_{i,1}(x_\alpha)=(-1)^{\alpha(h_i)}q_i^{\alpha(h_i)}L'_{i,1}(x_\alpha)$, we deduce that $L'_{i,1}(x_\alpha)\in U^+$ and therefore $L'_{i,1}(x)\in U^+$. Thus, the following four subspaces of $U^+$ coincide
$$U^+[i]; \ \text{ker}\delta^i;\ \{x\in U^+\mid L''_{i,1}(x)\in U^+\};\ \text{and}\ \{x\in U^+\mid L'_{i,1}(x)\in U^+\}.$$
By using $*$ , we have the following equal subspaces
$$U^+_*[i]; \ \text{ker}\delta_i;\ \{x\in U^+\mid L''_{i,-1}(x)\in U^+\};\ \text{and}\ \{x\in U^+\mid L'_{i,-1}(x)\in U^+\},$$
and by using $\varpi$, we obtain
$$U^-[i]=\{x\in U^-\mid L''_{i,1}(x)\in U^-\}=\{x\in U^-\mid L'_{i,1}(x)\in U^-\},$$
$$U_*^-[i]=\{x\in U^-\mid L''_{i,-1}(x)\in U^-\}=\{x\in U^-\mid L'_{i,-1}(x)\in U^-\}.$$

\vskip 3mm

Let $f_m=f(i,(j,l),m)$ and $f'_m=f'(i,(j,l),m)$. By direct calculation, we can get the following equations in $U^+$ easily as in \cite[Lemma 38.1.7]{Lusztig},
\begin{equation}\label{fm}
\begin{aligned}
& \rho(f_m)=1\otimes f_m+\sum_{t=0}^m  \prod_{h=0}^{m-t-1}(1-q_i^{2m-2h-2l\beta-2})q_i^{t(m-t)} f_t\otimes \mathtt a_i^{(m-t)},\\
& \rho(f'_m)=f'_m\otimes 1+\sum_{t=0}^m  \prod_{h=0}^{m-t-1}(1-q_i^{2m-2h-2l\beta-2})q_i^{t(m-t)}\mathtt a_i^{(m-t)}\otimes f'_t.
\end{aligned}
\end{equation}
Hence we have $\delta^i(f_m)=0$, $\delta_i(f_m)=(1-q_i^{2m-2l\beta-2})q_i^{m-1}f_{m-1}$ and
\begin{equation}\label{20}
\delta^{j,l}(f_m)= \prod_{h=0}^{m-1}(1-q_i^{2m-2h-2l\beta-2})\mathtt a_i^{(m)}=\prod_{h=1}^{m}(1-q_i^{-2(l\beta+1-h)})\mathtt a_i^{(m)}.
\end{equation}
Since $*\delta^{k,s}=\delta_{k,s}*$ for any $(k,s)\in I^\infty$, so we obtain
 $$\delta_i(f'_m)=0, \ \delta^i(f'_m)=(1-q_i^{2m-2l\beta-2})q_i^{m-1}f'_{m-1},$$
and
$$\delta_{j,l}(f'_m)=\prod_{h=1}^{m}(1-q_i^{-2(l\beta+1-h)})\mathtt a_i^{(m)}. $$

\vskip 3mm

\begin{lemma}\label{=}{\rm
Let $f_m=f(i,(j,l),m)$ and $g_m=g(i,(j,l),m)$, we have $$\{f_m,g_m\}=\{L''_{i,1}f_m,L''_{i,1}g_m\}.$$
\begin{proof}
Note that $\{ \mathtt B_iU^-,f_m\}=0$, so we have by (\ref{20}) that
$$\begin{aligned}
&\{g_m,f_m\}=(-1)^mq_i^{m(l\beta-m+1)}\{\mathtt b_{jl}\mathtt B_i^{(m)},f_m\}\\
&\phantom{\{g_m,f_m\}}=(-1)^mq_i^{m(l\beta-m+1)}\tau_{jl}\prod_{h=1}^{m}(1-q_i^{-2(l\beta+1-h)})\{\mathtt B_i^{(m)},\mathtt a_i^{(m)}\}.
\end{aligned}$$
Since $\{\mathtt a_i^{(m)},\mathtt a_i^{(m)}\}=\tau_i^mq_i^{m(m-1)/2}[m]_i^{-1}$, we have
$$\{\mathtt B_i^{(m)},\mathtt a_i^{(m)}\}=q_i^{m(m-1)/2}(q_i^{-1}-q_i)^{-m}{[m]_i}^{-1}.$$
Hence
$$\begin{aligned}& \{g_m,f_m\}=\tau_{jl}{\begin{bmatrix} l\beta \\ m \end{bmatrix}}_i=\{g_{l\beta-m},f_{l\beta-m}\}=\{*g_{l\beta-m},*f_{l\beta-m}\}\\
&\phantom{\{g_m,f_m\}}=\{g'_{l\beta-m},f'_{l\beta-m}\}=\{L''_{i,1}g_m,L''_{i,1}f_m\}.
\end{aligned}$$
The lemma is proved.
\end{proof}
}\end{lemma}

\vskip 3mm

Let $n>0$ and $i\in I^{\text{re}}$, we define the linear maps $\delta_{ni},\delta^{ni}:U^{+}\rightarrow U^{+}$ by
$$\begin{aligned}
& \rho(x)=x\otimes 1+\sum_{n>0}\delta_{ni}(x)\otimes \mathtt a_i^{(n)}+\text{terms of bidegree not in}\ Q_+\times\N\alpha_i,\\
& \rho(x)=1\otimes x+\sum_{n>0}\mathtt a_i^{(n)}\otimes \delta^{ni}(x)+\text{terms of bidegree not in}\ \N\alpha_i\times Q_+,
\end{aligned}$$
where $x$ is a homogeneous element of $U^+$. Let $|x|=\mu$ and let $(j,l)\in I^\infty$ with $i\neq (j,l)$, we have
$$\begin{aligned}
& \delta_{ni}(\mathtt a_ix)=\mathtt a_i\delta_{ni}(x)+[n]_iq^{(\alpha_i,\mu-(n-1)\alpha_i)}\delta_{(n-1)i}(x),\\
& \delta^{ni}(x\mathtt a_i)=\delta^{ni}(x)\mathtt a_i+[n]_iq^{(\alpha_i,\mu-(n-1)\alpha_i)}\delta^{(n-1)i}(x),
\end{aligned}$$
and $$\delta_{ni}(\mathtt a_{jl}x)=\mathtt a_{jl}\delta_{ni}(x),\ \delta^{ni}(x\mathtt a_{jl})=\delta^{ni}(x)\mathtt a_{jl}.$$

Note that, we have by induction on $n$ that
$$\begin{aligned}
& (\delta_i)^n(\mathtt a_ix)=\mathtt a_i(\delta_i)^n(x)+[n]_iq^{(\alpha_i,\mu)}q_i^{-(n-1)}(\delta_i)^{n-1}(x),\\
& (\delta^{i})^n(x\mathtt a_i)=(\delta^{i})^n(x)\mathtt a_i+[n]_iq^{(\alpha_i,\mu)}q_i^{-(n-1)}(\delta^i)^{n-1}(x),
\end{aligned}$$
and $$(\delta_{i})^n(\mathtt a_{jl}x)=\mathtt a_{jl}(\delta_i)^n(x),\ (\delta^i)^n(x\mathtt a_{jl})=(\delta^{i})^n(x)\mathtt a_{jl}.$$
It follows that
\begin{equation}\label{21}
(\delta_i)^n=q_i^{n(n-1)/2}\delta_{ni} \ \text{and} \ (\delta^i)^n=q_i^{n(n-1)/2}\delta^{ni}.
\end{equation}

\vskip 3mm

Fix $i\in I^{\text{re}}$, let $(j,l)\in I^\infty$ with $(j,l)\neq i$. Note that, the following set of elements are linearly independent in $U^+$
$$\{\mathtt a_i^{(n)}; \mathtt a_i^{(r)}\mathtt a_{jl}\mathtt a_i^{(s)}\mid n,r,s\geq 0 \ \text{and}\ r+s\leq -la_{ij}\}.$$
It can be extended to be a basis of $U^+$, by adding a set $\{u_w\}_{w\in S}$ of monomials in $U^+$.

Hence for any $x\in U^+$, we can write $\rho(x)$ into
\begin{equation}\label{ex}
\rho(x)=1\otimes x+\sum_{n>0}\mathtt a_i^{(n)}\otimes \delta^{ni}(x)+\sum_{r=0}^{l\beta}\sum_{s=0}^{l\beta-r}\mathtt a_i^{(r)}\mathtt a_{jl}\mathtt a_i^{(s)}\otimes \delta^{ri;(j,l);si}(x)+\sum_{w\in S}u_{w}\otimes x_{w}.
\end{equation}
Here $\beta=-a_{ij}$ as usual.

\vskip 3mm

\begin{lemma}\label{mi}{\rm
Let $x\in U^+$ and $y\in U^-$. For any $m\in \Z$, we have
\begin{itemize}
\item [(i)] $\{g_my,x\}=\{g_m,f_m\}\{y,\delta^{(j,l);mi}(x)\}.$
\item [(ii)] $\{g'_my,x\}=\{g'_m,f'_m\}\{y,\delta^{mi;(j,l)}(x)\}.$
\end{itemize}
\begin{proof}
We shall prove (i) first. Since $f_m=g_m=0$ when $m> l\beta$ or $m<0$, we may assume $0\leq m\leq l\beta$. Let $u$ be a monomial in $U^+_{m\alpha_i+l\alpha_j}$. If $u$ contains no $\mathtt a_{jl}$, then we can write $u=\mathtt a_i^{t}\mathtt a_{jk}u'$ for some $t\geq 0$, $k<l$ and a monomial $u'$ of $U^+_{(m-t)\alpha_i+(l-k)\alpha_j}$. Note that
$$
 \rho(\mathtt a_i^{(r)}\mathtt a_{jl}\mathtt a_i^{(m-r)})=(\sum_{a+b=r}q_i^{ab}\mathtt a_i^{(a)}\otimes \mathtt a_i^{(b)})(\mathtt a_{jl}\otimes 1+1\otimes \mathtt a_{jl})(\sum_{a'+b'=m-r}q_i^{a'b'}\mathtt a_i^{(a')}\otimes \mathtt a_i^{(b')}),$$
and the right hand side is in
$$ \sum_{\substack{ a+b=r\\ a'+b'=m-r}}\Q(q)\mathtt a_i^{(a)}\mathtt a_{jl}\mathtt a_i^{(a')}\otimes \mathtt a_i^{b+b'}+\sum_{\substack{ a+b=r\\ a'+b'=m-r}}\Q(q)\mathtt a_i^{a+a'}\otimes \mathtt a_i^{(b)}\mathtt a_{jl}\mathtt a_i^{(b')}.
$$
Hence we obtain $$\{\mathtt a_i^{(r)}\mathtt a_{jl}\mathtt a_i^{(m-r)},u\}=\{ \rho(\mathtt a_i^{(r)}\mathtt a_{jl}\mathtt a_i^{(m-r)}),\mathtt a_i^{t}\mathtt a_{jk}\otimes u' \}=0,$$
and therefore $\{g_m,u\}=0$ for such $u$.

Since $\{g_m,\mathtt a_iU^+\}=0$, we get $\{g_m,f_m\}=\{g_m,\mathtt a_{jl}\mathtt a_i^{(m)}\}$. On the other hand, by the above argument, we have $\{g_m,u_w\}=0$ for any $w\in S$, which implies
$$\{g_my,x\}=\{g_m\otimes y,\mathtt a_{jl}\mathtt a_i^{(m)}\otimes \delta^{(j,l);mi}(x)\}=\{g_m,f_m\}\{y,\delta^{(j,l);mi}(x)\}.$$ The proof of (ii) is similar.
\end{proof}
}\end{lemma}

\vskip 3mm

Let $x\in U^+_{\mu}$ and $x'\in U^+_{\mu'}$. Let $0\leq n\leq l\beta$. The expression (\ref{ex}) yields
\begin{equation}
\begin{aligned}
& \delta^{(j,l);ni}(xx')=\delta^{(j,l);ni}(x)x'+q^{(\mu,l\alpha_j+n\alpha_i)}x\delta^{(j,l);ni}(x')\\
& \phantom{\delta}+\sum_{t=0}^{n-1}q^{(\mu-l\alpha_j-t\alpha_i,(n-t)\alpha_i)}{\begin{bmatrix}n\\ t\end{bmatrix}}_i\delta^{(j,l);ti}(x)\delta^{(n-t)i}(x'),
\end{aligned}
\end{equation}
and
\begin{equation}\label{z2}
\begin{aligned}
& \delta^{ni;(j,l)}(xx')=\delta^{ni;(j,l)}(x)x'+q^{(\mu,n\alpha_i+l\alpha_j)}x\delta^{ni;(j,l)}(x')\\
& \phantom{\delta}+\sum_{t=0}^{n-1}q^{(\mu-(n-t)\alpha_i,t\alpha_i+l\alpha_j)}{\begin{bmatrix}n\\ t\end{bmatrix}}_i\delta^{(n-t)i}(x)\delta^{ti;(j,l)}(x').
\end{aligned}
\end{equation}
In particular, if $x'\in U^+[i]$, we have $\delta^{ni}(x')=0$ for all $n>0$, which implies
\begin{equation}\label{sp}
\delta^{(j,l);ni}(xx')=\delta^{(j,l);ni}(x)x'+q^{(\mu,l\alpha_j+n\alpha_i)}x\delta^{(j,l);ni}(x').
\end{equation}

\vskip 3mm
Let $f_m=f(i,(j,l),m)$ and $f'_m=f'(i,(j,l),m)$. By ($\ref{fm}$), we have
$$ \delta^{(j,l);ni}(f_m)=\begin{cases}\gamma_{mn}\mathtt a_i^{(m-n)}& \text{if}\ n\leq m,\\ 0 & \text{if} \ n>m,\end{cases} \  \text{and} \ \ \delta^{ni;(j,l)}(f'_m)=\begin{cases}1\ \ \text{if}\ n=m,\\ 0 \ \ \text{otherwise},\end{cases}$$
where $$\gamma_{mn}=\prod_{h=0}^{m-n-1}(1-q_i^{2m-2h-2l\beta-2})q_i^{n(m-n)}.$$

\vskip 3mm

Denote by $P_i:U^+\rightarrow U^+[i]$ the projection obeys the decomposition $U^+=\mathtt a_iU^+\oplus U^+[i]$ given in Lemma \ref{decom}. We have by the definition
$$P_i(xx')=P_i(P_i(x)x') \ \text{for any}\ x,x'\in U^+.$$
In particular, we have a more simpler formula $P_i(xx')=P_i(x)x'$ when $x'\in U^+[i]$.

\vskip 3mm

\begin{lemma}\label{rx}{\rm
Let $x\in U^+[i]$ with $|x|=\mu$ and let $n>0$, we have
$$P_i(x\mathtt a_i^n)=\frac{q^{n(\mu,\alpha_i)}q_i^{n(n+1)}}{(q_i-q_i^{-1})^n}L'_{i,-1}(\delta^{i})^nL''_{i,1}(x).$$
\begin{proof}
We use induction on $n$. Assume that $n=1$. Since $L''_{i,1}(x)\in U^+_*[i]$, we have $\delta_i(L''_{i,1}(x))=0$ and so
$$L''_{i,1}(x)\mathtt B_i-\mathtt B_iL''_{i.1}(x)=-K_i^{-1}\delta^i(L''_{i,1}(x))/ (q_i-q_i^{-1}).$$

Applying $L'_{i,-1}$ to both sides, we get
$$-x\mathtt a_iK_i+\mathtt a_iK_ix=-K_iL'_{i,-1}\delta^iL''_{i,1}(x)/(q_i-q_i^{-1}),$$
and by applying $K^{-1}$, we obtain $$-q^{-(\mu+\alpha_i,\alpha_i)}x\mathtt a_i+q^{-(\alpha_i,\alpha_i)}\mathtt a_ix=-L'_{i,-1}\delta^iL''_{i,1}(x)/(q_i-q_i^{-1}).$$
It follows that $$x\mathtt a_i=q^{(\mu,\alpha_i)}\mathtt a_ix+\frac{q^{(\mu+\alpha_i,\alpha_i)}}{q_i-q_i^{-1}}L'_{i,-1}\delta^iL''_{i,1}(x).$$
Note that, according to the second formula in (\ref{fm}), we have $\rho(U^+_*[i])\subseteq U^+\otimes U^+_*[i]$, which implies $\delta^iL''_{i,1}(x)\in U^+_*[i]$. Hence
$$P_i(x\mathtt a_i)=\frac{q^{(\mu+\alpha_i,\alpha_i)}}{q_i-q_i^{-1}}L'_{i,-1}\delta^iL''_{i,1}(x).$$
This proves our assertion for $n=1$. For the induction step, assume the lemma is true for $n$, then we have
$$\begin{aligned}
& P_i(x\mathtt a_i^{n+1})=P_i(P_i(x\mathtt a_i^n)\mathtt a_i)=\frac{q^{(\mu+(n+1)\alpha_i,\alpha_i)}}{q_i-q_i^{-1}}L'_{i,-1}\delta^iL''_{i,1}(P_i(x\mathtt a_i^n))\\
&\phantom{P_i(x\mathtt a_i^{n+1})}=\frac{q^{(n+1)(\mu,\alpha_i)}q_i^{(n+1)(n+2)}}{(q_i-q_i^{-1})^{n+1}}L'_{i,-1}(\delta^i)^{n+1}L''_{i,1}(x)
\end{aligned}$$
as desired.
\end{proof}
}\end{lemma}

\vskip 3mm

Now, we shall prove our main theorem by using the method given in \cite[Chapter 8A]{Jantzen}. Note that there is no special restriction on  the values of $\tau_i$'s for $i\in I^{\text{re}}$, but only ask them to take values  in $1+q^{-1}\Z_{\geq0}[[q^{-1}]]$, hence the argument in \cite[Lemma 38.2.1]{Lusztig} is not available for our case.

\vskip 3mm

\begin{theorem}\label{remain}{\rm
For any $x\in U^+[i], y\in U^-[i]$, we have
$$\{L''_{i,1}(x),L''_{i,1}(y)\}=\{x,y\}.$$
\begin{proof}
Assume that our assertion holds for a given $y$ in $U^-[i]$ and arbitrary $x$ in $U^+[i]$. By Lemma \ref{mi}(ii), we have
$$\begin{aligned}
& \{L''_{i,1}(g_my),L''_{i,1}(x)\}=\{g'_{l\beta-m}L''_{i,1}(y),L''_{i,1}(x)\}\\
& =\{g'_{l\beta-m}, f'_{l\beta-m}\}\{L''_{i,1}(y),\delta^{(l\beta-m)i;(j,l)}L''_{i,1}(x)\}\\
& =\{g'_{l\beta-m}, f'_{l\beta-m}\}\{y,L'_{i,-1}\delta^{(l\beta-m)i;(j,l)}L''_{i,1}(x)\},
\end{aligned}$$
where the last equality follows from the fact that $\delta^{(l\beta-m)i;(j,l)}L''_{i,1}(x)\in U^+_*[i]$. On the other hand, by Lemma \ref{=} and Lemma \ref{mi}(i), we have $$\{g_my,x\}=\{g_m,f_m\}\{y,\delta^{(j,l);mi}(x)\}=\{g'_{l\beta-m}, f'_{l\beta-m}\}\{y,\delta^{(j,l);mi}(x)\}.$$

Since $\{y,\mathtt a_iU^+\}=0$. To show $\{L''_{i,1}(g_my),L''_{i,1}(x)\}=\{g_my,x\}$, it is enough to show
$$L'_{i,-1}\delta^{(l\beta-m)i;(j,l)}L''_{i,1}(x)\equiv \delta^{(j,l);mi}(x) \ \text{mod}\ \mathtt a_iU^+$$
for all $x\in U^+[i]$. This is equivalent to
\begin{equation}\label{pi}
\delta^{(l\beta-m)i;(j,l)}L''_{i,1}(x)=L''_{i,1}P_i \delta^{(j,l);mi}(x).
\end{equation}
for all $x\in U^+[i]$.

The identity (\ref{pi}) holds for $x=1$, since both sides equal $0$ in this case. Assume (\ref{pi}) holds for $x\in U^+[i]$ with $|x|=\mu$, and let $x'=f(i,(j',l'),s)$. If $(j',l')\neq (j,l)$, then the first equality in (\ref{fm}) yields
$\delta^{(j,l);mi}(x')=0$. By (\ref{sp}), we have $\delta^{(j,l);mi}(xx')=\delta^{(j,l);mi}(x)x'$,
which implies
$$P_i\delta^{(j,l);mi}(xx')=P_i(\delta^{(j,l);mi}(x))x',$$
and so $$L''_{i,1}P_i\delta^{(j,l);mi}(xx')=L''_{i,1}P_i\delta^{(j,l);mi}(x)L''_{i,1}(x').$$
On the other hand, by (\ref{z2}), we have $$\delta^{(l\beta-m)i;(j,l)}L''_{i,1}(xx')=\delta^{(l\beta-m)i;(j,l)}L''_{i,1}(x)L''_{i,1}(x'),$$
hence (\ref{pi}) holds for $xx'$ in this case.

If $(j',l')=(j,l)$, then $x'=f_s=f(i,(j,l),s)$. By using (\ref{sp}), we obtain
$$\delta^{(j,l);mi}(xf_s)=\delta^{(j,l);mi}(x)f_s+X,$$
where $X=0$ when $m>s$, and $X=q^{(\mu,l\alpha_j+m\alpha_i)}\gamma_{sm}x\mathtt a_i^{(s-m)}$ when $m\leq s$.

By (\ref{z2}), we have
$$\delta^{(l\beta-m)i;(j,l)}L''_{i,1}(xf_s)=\delta^{(l\beta-m)i;(j,l)}L''_{i,1}(x)f'_{l\beta-s}+Y,$$
where $Y=0$ when $m>s$, and when $m\leq s$,
$$Y=q^{(r_i\mu-(s-m)\alpha_i,(l\beta-s)\alpha_i+l\alpha_j)}{\left[ \substack{ l\beta-m\\l\beta-s}\right]}_i\delta^{(s-m)i}L''_{i,1}(x).$$ Here we understand $\delta^{0i}(L''_{i,1}(x))=L''_{i,1}(x)$.
According to our assumption, we have $$L''_{i,1}\left(P_i\delta^{(j,l);mi}(x)f_s\right)=\delta^{(l\beta-m)i;(j,l)}L''_{i,1}(x)f'_{l\beta-s}.$$

Hence in order to show
$$\delta^{(l\beta-m)i;(j,l)}L''_{i,1}(xf_s)=L''_{i,1}P_i\delta^{(j,l);mi}(xf_s),$$
it suffices to prove that $L''_{i,1}P_i(X)=Y$ in the case of $m\leq s$. In this case, using the equality given Lemma \ref{rx} and in (\ref{21}), we get
$$\begin{aligned}
& L''_{i,1}P_i(X)=q^{(\mu,l\alpha_j+m\alpha_i)}\gamma_{sm}L''_{i,1}P_i(x\mathtt a_i^{(s-m)})\\
& \ \ =\frac{1}{[s-m]_i!}\frac{q^{(s-m)(\mu,\alpha_i)}q^{(\mu,l\alpha_j+m\alpha_i)}}{(q_i-q_i^{-1})^{(s-m)}}q_i^{(s-m)(3(s-m)+1)/2}\gamma_{sm}\delta^{(s-m)i}L''_{i,1}(x).
\end{aligned}$$
Note that, we have
$$ \gamma_{sm}
=q_i^{m(s-m)+(s-l\beta)(s-m)-(s-m)(s-m+1)/2}\prod_{h=1}^{s-m}(q_i^{l\beta-s+h}-q_i^{-l\beta+s-h}).
$$
Thus, $$ L''_{i,1}P_i(X)=q^{(\mu,s\alpha_i+l\alpha_j)}q_i^{(2s-l\beta)(s-m)}{\left[ \substack{ l\beta-m\\l\beta-s}\right]}_i\delta^{(s-m)i}L''_{i,1}(x)=Y.$$
The theorem is proved.
\end{proof}
}\end{theorem}

\vskip 3mm
Let $P'_i:U^+\rightarrow U^+_*[i]$ be the projection alongs the decomposition $U^+= U^+_*[i]\oplus U^+\mathtt a_i$. We have the following result as a direct corollary of the previous theorem.
\vskip 3mm

\begin{corollary}{\rm
Let $x\in U^+[i]$, we have
\begin{equation}\label{go}
(P'_i\otimes \text{id})\circ\rho \circ L''_{i,1}(x)=(L''_{i,1}\otimes L''_{i,1})\circ (\text{id}\otimes P_i)\circ \rho(x).
\end{equation}
\begin{proof} Note that both sides in the equation (\ref{go}) belong to $U^+_*[i]\otimes U^+_*[i]$. Let $y,y'\in U^-[i]$, we have by the definition
$$\begin{aligned}& \{yy',x\}=\{y\otimes y',\rho(x)\}=\{y\otimes y',(\text{id}\otimes P_i)\circ\rho(x)\}\\
&\phantom{\{yy',x\}}=\{L''_{i,1}(y)\otimes L''_{i,1}(y'),(L''_{i,1}\otimes L''_{i,1})\circ (\text{id}\otimes P_i)\circ\rho(x)\},\end{aligned}$$
where the last equality follows from Theorem \ref{remain}.
On the other hand, since $yy'\in U^-[i]$,  we obtain by using Theorem \ref{remain} directly
$$\begin{aligned}
& \{yy',x\}=\{L''_{i,1}(y)L''_{i,1}(y'),L''_{i,1}(x)\}=\{L''_{i,1}(y)\otimes L''_{i,1}(y'),\rho L''_{i,1}(x)\}\\
& \phantom{\{yy',x\}}=\{L''_{i,1}(y)\otimes L''_{i,1}(y'),(P'_i\otimes \text{id})\circ\rho \circ L''_{i,1}(x)\}.
\end{aligned}$$

Let $u=(P'_i\otimes \text{id})\circ\rho \circ L''_{i,1}(x)-(L''_{i,1}\otimes L''_{i,1})\circ (\text{id}\otimes P_i)\circ \rho(x)$. The above argument implies that
$$\{U^-_*[i]\otimes U^-_*[i],u\}=0.$$

Since $U^-=U^-_*[i]\oplus U^-\mathtt B_i$ and $\{U^-\mathtt B_i,u\}=0$, thus we deduce that $\{U^-\otimes U^-,u\}=0$.  Note that the bilinear form $\{ \ , \ \}$ is non-degenerated on $U^+\otimes U^+$, hence $u=0$ and our assertion follows.
\end{proof}
}\end{corollary}

\end{document}